\documentclass[12pt,leqno]{article}
\usepackage{amssymb}
\usepackage{amscd}
\usepackage{amsmath}

\newcommand{\Q}{\mathbb Q}

\newcommand{\R}{\mathbb R}
\newcommand{\C}{\mathbb C}

\newcommand{\A}{\mathbb A}

\makeatletter
\def\l@section{\@tocline{1}{4pt}{1pc}{}{}}
\def\l@subsection{\@tocline{2}{0pt}{2pc}{5pc}{}}
\makeatother

\title{Existence of Ramanujan primes for GL$(3)$}

\author{Dinakar Ramakrishnan \\253-37 Caltech,
Pasadena, CA 91125}
\date{}
\begin{document}
\maketitle \centerline{To Joe Shalika with admiration}
\bigskip

\section*
{\bf Introduction}

\medskip

Let $\pi$ be a cusp form on GL$(n)/\Q$, i.e., a cuspidal
automophic representation of GL$(n, \A)$, where $\A$ denotes the
adele ring of $\Q$. We will say that a prime $p$ is a {\bf
Ramanujan prime for ${\bf \pi}$} iff the corresponding $\pi_p$ is
{\it tempered}. The local component $\pi_p$ will necessarily be
{\it unramified} for almost all $p$, determined by an unordered
$n$-tuple $\{\alpha_{1,p}, \alpha_{2,p}, \dots, \alpha_{n,p}\}$ of
non-zero complex numbers, often represented by the corresponding
diagonal matrix $A_p(\pi)$ in GL$(n, \C)$, unique up to
permutation of the diagonal entries. The $L$-factor of $\pi$ at
$p$ is given by
$$
L(s, \pi_p) \, = \, {\rm det}(I - A_p(\pi)p^{-s})^{-1} \, = \,
\prod\limits_{j=1}^n (1-\alpha_{j,p}p^{-s})^{-1}.
$$
As $\pi$ is unitary, $\pi_p$ is tempered (in the unramified case)
iff each $\alpha_{j,p}$ is of absolute value $1$. It was shown in
[Ra] that for ${\bf n=2}$, the set $\mathcal R(\pi)$ of Ramanujan
primes for $\pi$ is of lower density at least $9/10$. When one
applies in addition the deep recent results of H.~Kim and
F.~Shahidi ([KSh]) on the symmetric cube and the symmetric fourth
power liftings for GL$(2)$, the lower bound improves from $9/10$
to $34/35$ ({\it loc. cit.}), which is $0.971428\dots$.

\medskip

For $n > 2$ there is a dearth of results for general $\pi$, though
for cusp forms of {\it regular infinity type}, assumed for $n
> 3$ to have a discrete series component at a finite place,
one knows by [Pic] for $n=3$, which relies on the works of
J.~Rogawski, et al, and by the work of Clozel ([C$\ell$]) for $n >
3$ (see also the non-trivial refinement due to Harris and Taylor
([HaT])), that {\it all} the unramified primes are Ramanujan
primes for $\pi$. (Of course for $n =2$, a regular cusp form is
necessarily holomorphic of weight $k \geq 2$, and it is known, by
Eichler-Shimura for $k=2$ and by Deligne for $k > 2$, that every
prime is a Ramanujan prime in this case; ditto for $k=1$ by the
work of Deligne and Serre.) One is interested in knowing whether
there exists even one Ramanujan prime for general $\pi$ on
GL$(n)$. Thanks to the work of Kim and Shahidi, one sees the
importance of knowing a positive answer to such a question.

\medskip

The main result of this Note is the following:

\medskip

\noindent{\bf Theorem A} \, \it Let $\pi$ be a cusp form on
GL$(3)/\Q$. Then there exist infinitely many Ramanujan primes for
$\pi$.
\rm

\medskip

Let us now explain the main issues behind its proof. One can show
(see section 1) that at any prime $p$ where $\pi$ is unramified,
{\it if} the coefficient $a_p(\pi)$ is bounded in absolute value
by $1$, {\it then} $\pi_p$ is tempered. A general result proved in
[Ra] for GL$(n)$ implies that for any real number $b
> 1$, the set of $p$ where $\vert a_p(\pi)\vert \leq b$ is
infinite, even of lower Dirichlet density $\geq 1 -
\frac{1}{b^2}$. But this gives us nothing for $b =1$. Our aim here
is to show that for infinitely many primes $p$, $a_p(\pi)$ is
indeed bounded in absolute value by $1$. The key idea is to
exploit the adjoint $L$-function (whose definition makes sense for
$\pi$ on GL$(n)$ for any $n$):
$$
L(s, \pi; Ad) \, =  \, \frac{L(s, \pi \times \overline
\pi)}{\zeta(s)},
\leqno(0.1)
$$
where $L(s, \pi \times \overline \pi)$ is the Rankin-Selberg
$L$-function of the pair $(\pi, \overline \pi)$. (As usual,
$\overline \pi$ signifies the complex conjugate of $\pi$, which,
by the unitarity of $\pi$, is the same as the contragredient of
$\pi$.) One knows (see [HRa], Lemma $a$ of section 2) that $L(s,
\pi \times \overline \pi)$ is of {\it positive type}, i.e., the
Dirichlet series defined by its logarithm has non-negative
coefficients. The proof of Theorem A relies on the following

\medskip

\noindent{\bf Proposition B} \, \it Let $\pi$ be a cusp form on
GL$(n)/\Q$. Then for any finite set $S$ of primes containing
infinity, the incomplete adjoint $L$-function $L^S(s, \pi, Ad)$ is
not of positive type. \rm

\medskip

The proof given here of this Proposition, and hence of Theorem A,
will work over any number field $F$ having no real zeros in the
interval $(0,1)$. In the case of real zeros one has to proceed
differently.

\medskip

When I was at Hopkins as an Assistant Professor during 1983-85, I
learnt a lot from Joe Shalika about the $L$-functions of GL$(n)$.
It is a pleasure to dedicate this article to him. Thanks are due
to Jeff Lagarias and Freydoon Shahidi for making comments on an
earlier version of this article, which led to an improvement of
the exposition, and to the NSF for financial support through the
grant DMS--0100372.

\vskip 0.2in

\section{Why Proposition B implies Theorem A}

\bigskip

Let $\pi \, \simeq \, \pi_\infty \otimes(\otimes_p' \pi_p)$ be a
cuspidal automorphic representation of GL$(3, \A) = {\rm GL}(3,
\R) \times {\rm GL}(3,\A_f)$. Let $S_0$ be the set of primes $p$
where $\pi_p$ is unramified and tempered, and let $S_1$ be the
finite set of primes where $\pi_p$ is ramified. Put
$$
S \, = \, S_0 \cup S_1 \cup \{\infty\}. \leqno(1.1)
$$

For any $L$-function with an Euler product $\prod_v L_v(s)$ over
$\Q$, put
$$
L^S(s) \, = \, \prod_{p \notin S} L_p(s),
\leqno(1.2)
$$
which we call the {\it incomplete Euler product} relative to, or
outside, $S$.

\medskip

Pick any $p$ outside $S$ and consider the Langlands class
$$
A_p \, = \, A_p(\pi) \, = \, \{\alpha_{1,p}, \alpha_{2,p},
\alpha_{3,p}\}. \leqno(1.3)
$$
As $\pi_p$ is by assumption non-tempered, there is a non-zero real
number $t$ and a complex number $u$ of absolute value $1$ such
that, after possibly renumbering the $\alpha_{j,p}$,
$$
\alpha_{1,p} \, = \, up^t.
$$
On the other hand, by the unitarity of $\pi_p$, we must have
$$
\{\overline \alpha_{1,p}, \overline\alpha_{2,p},
\overline\alpha_{3,p}\} \, = \, \{\alpha_{1,p}^{-1},
\alpha_{2,p}^{-1}, \alpha_{3,p}^{-1}\}.
$$
This then implies that
$$
A_p \, = \, \{up^t, up^{-t}, w\}, \leqno(1.4)
$$
for some complex number $w$ of absolute value $1$. We may, and we
will, assume that $t$ is positive. Put
$$
u^{-1}w \, = \, e^{i\theta},
\leqno(1.6)
$$
for some $\theta \in [0, 2\pi) \subset \R$.

So we have
$$
\vert a_p \vert^2 \, = \, (p^t + p^{-t} + \cos\theta)^2 +
\sin^2\theta \, = \, 3 + p^{2t} + p^{-2t} + 2\cos\theta(p^t +
p^{-t}). \leqno(1.6)
$$

\medskip

Now let us look at the adjoint $L$-function. By definition,
$$
L^S(s, \pi; Ad) \, = \, \frac{\prod_{p \notin S}L(s, \pi_p
\overline \pi_p)}{\prod_{p \notin S} (1-p^{-s})^{-1}}.
\leqno(1.7)
$$
So for any $p$ outside $S$, the Langlands class of the Adjoint
$L$-function is
$$
A_p(\pi; Ad) \, = \, A_p \otimes \overline A_p - \{1\}.
$$
Applying (1.4) and (1.5), we get
$$
A_p(\pi; Ad) \, = \, \{p^{2t}, p^{-2t}, 1, 1, u\overline wp^t,
u\overline w p^{-t}, \overline uwp^t, \overline uwp^{-t}\}.
\leqno(1.8)
$$
and
$$
a_p(\pi; Ad) \, = \, {\rm tr}(A_p(\pi; Ad)) \, = \,
2+p^{2t}+p^{-2t}+ 2\cos\theta(p^t+p^{-t}).
$$
Consequently,
$$
\log L^S(s, \pi; Ad) \, = \, \sum\limits_{p \notin S}
\sum\limits_{m \geq 1} \frac{a_{p^m}(\pi; Ad)}{p^{ms}},
\leqno(1.9)
$$
where (by (1.8) and (1.6))
$$
a_{p^m}(\pi; Ad) \, = \, 2+p^{2mt}+p^{-2mt}+ 2\cos
m\theta(p^{mt}+p^{-mt}).
\leqno(1.10).
$$
Since
$$
p^{mt}+p^{-mt} \, \geq \, 2,
$$
and since
$$
a_{p^m}(\pi; Ad) \, = \, (p^{mt}+p^{-mt})((p^{mt}+p^{-mt})+2\cos
m\theta)),
$$
we get the following

\medskip

\noindent{\bf Lemma 1.11} \, \it Let $\pi$ be a cusp form on
GL$(3)/\Q$ and $S$ the set of primes containing $\infty$, the
primes where $\pi$ is ramified and the Ramanujan primes for $\pi$.
Then $L^S(s, \pi; Ad)$ is of positive type. \rm

\medskip

But if $S_0$, and hence $S$, is finite, this Lemma contradicts the
conclusion of Proposition B. Hence the set of Ramanujan primes for
$\pi$ must be infinite, once we accept Proposition B.

\vskip 0.2in

\section{Proof of Proposition B}

\medskip

In this section $\pi$ will be a unitary, cuspidal automorphic
representation of GL$(n, \A)$. At each place $v$, the local factor
of $L(s, \pi; Ad)$ is given by
$$
L(s, \pi_v; Ad) \, = \, \frac{L(s, \pi_v \times \overline
\pi_v)}{\zeta_v(s)},
\leqno(2.1)
$$
where $\zeta_v(s)$ is $(1-p^{-s})^{-1}$ if $v$ is a finite place
defined by a prime $p$, and it equals $\pi^{s/2}\Gamma(s/2)$ if
$v$ is the archimedean place. By convention, $\zeta(s)$ is the
product of $\zeta_v(s)$ over all the {\it finite} $v$, while all
the other automorphic $L$-functions occurring in this paper will
also involve the archimedean factors.

The $L$-group of GL$(n)$ is GL$(n, \C)$, and the Euler factor
$L(s, \pi; Ad)$ is associated to the representation
$$
Ad: \, {\rm GL}(n, \C) \, \rightarrow \, {\rm GL}(n^2-1, \C),
\leqno(2.2)
$$
given by composing the natural projection of GL$(n, \C)$ onto
PGL$(n, \C)$ with the $(n^2-1)$-dimensional {\it Adjoint
representation} of PGL$(n, \C)$, whence the notation $Ad$. In any
case, we have for {\it every} $v$:
$$
L(s, \pi_v; Ad) \, \ne \, 0 \, \, \forall \, s \in \C.
\leqno(2.3)
$$
One way to see this will be to use the local Langlands
correspondence, established recently by Harris-Taylor ([HaT]) and
Henniart ([He]), associating to each $\pi_v$ an $n$-dimensional
representation $\sigma_v$ of $W_{F_v} \times {\rm SL}(2, \C)$
(resp. $W_{F_v}$) for $v$ finite (resp. infinite). (Here $W_{F_v}$
denotes as usual the Weil group of $F_v$.) Since this
correspondence preserves the local factors of pairs and matches
identifie the central character of $\pi_v$ with the determinant of
$\sigma_v$, one gets in particular,
$$
L(s, \pi_v; Ad) \, = \, L(s, Ad(\sigma_v)),
$$
where $Ad(\sigma_v)$ denotes $\sigma_v \otimes \sigma_v^\vee
\ominus 1$, which is a {\it genuine} representation because the
trivial representation occurs in $\sigma_v \otimes \sigma_v^\vee
\, \simeq \, $End$(\sigma_v)$. It is well known that for any
representation $\tau_v$ of $W_{F_v} \times {\rm SL}(2, \C)$, such
as $Ad(\sigma_v)$, the associated $L$-factor has no zeros.

\medskip

Now let $S$ be any finite set of primes containing $\infty$ and
the primes where $\pi$ is ramified. Put
$$
L^S(s, \pi; Ad) \, = \, \prod_{v \notin S} L(s, \pi_v; Ad).
\leqno(2.4)
$$

\medskip

Suppose $L^S(s, \pi; Ad)$ is of positive type. Then by definition,
its logarithm defines a Dirichlet series with positive
coefficients, absolutely convergent in a right half plane. By the
theory of Landau, this Dirichlet series converges on $(\beta,
\infty)$, where $\beta$ is the largest real number where $\log
L^S(s, \pi; Ad)$ diverges. But such a point of divergence must be
a pole, and not a zero, of $L^S(s, \pi; Ad)$ because its logarithm
is positive in $(\beta, \infty)$.

\medskip

\noindent{\bf Lemma 2.5} \, \it Let $\beta$ be the largest real
number such that $L^S(s, \pi; Ad)$ converges for all real $s >
\beta$. Then
$$
\beta \, < \, 0.
$$
\rm

\medskip

{\it Proof of Lemma 2.5}. \, By the standard properties of the
Rankin-Selberg $L$-functions ([JPSS], [JS], [Sh1-3], [MW] -- see
also [BRa]), $L^S(s, \pi \times \overline \pi)$ is invertible for
$\Re(s) > 1$ and admits a meromorphic continuation to the whole
$s$-plane with a unique {\it simple} pole at $s=1$. The same
properties hold of course for $\zeta^S(s)$; so $L^S(s, \pi; Ad)$
has no pole in $\Re(s) \geq 1$. In other words we have $\beta <
1$. Moreover, one knows that $\zeta^S(s)$ is {\bf non-zero} on
$(0, 1) \subset \R$. Hence we have
$$
\beta \, \leq \, 0.
$$
Now let us look at the point $s=0$. By definition,
$$
L^S(s, \pi; Ad) \, = \, \frac{L^\infty(s, \pi \times \overline
\pi)}{\zeta(s)\prod_{v \in S - \{\infty\}} L(s, \pi_v; Ad)}.
$$
The numerator on the right has no pole at $s=0$, and $\zeta(s)$
does not vanish at $s=0$. Moreover, as we noted above, the local
factors $L(s, \pi_v; Ad)$ have no zeros. Consequently, $L^S(s,
\pi; Ad)$ has no pole at $s=0$, and this proves the Lemma.

\bigskip

\noindent{\bf Lemma 2.6} \, \it $L(s, \pi_\infty; Ad)$ has a pole
at $s=0$. \rm

\medskip

{\it Proof of Lemma 2.6}. \, There exist complex numbers $z_1,
z_2, z_3$ such that
$$
L(s, \pi_\infty) \, = \, \prod\limits_{j=1}^3
\Gamma_\R(s+z_j+\delta_j),
\leqno(2.7)
$$
with $\delta_j \in \{\pm 1\}, \, \forall j$, and
$$
\Gamma_\R(s) \, = \, \pi^{-s/2}\Gamma(s/2).
$$
By the unitarity of $\pi_\infty$, we see that either all the $z_j$
have absolute value $1$, in which case $\pi_\infty$ is tempered,
or exactly one of the $z_j$, say $z_1$, has absolute value $1$,
and moreover,
$$
z_2 = u+t, \, z_3 = u-t,
$$
for some positive real number $t$ and a complex number $u$ of
absolute value $1$. In either case we see that the set
$$
B(\pi_\infty; Ad) : = \, \{z_1, z_2, z_3\} \cup \{\overline z_1,
\overline z_2, \overline z_3\} - \{0\}
\leqno(2.8)
$$
contains $0$. The standard yoga of Langlands $L$-functions
furnishes the identity
$$
L(s, \pi_\infty; Ad) \, = \, \prod\limits_{z \in B(\pi_\infty;
Ad)} \, \Gamma_\R(s+z).
\leqno(2.9)
$$
Recall that $\Gamma_\R(s)$ never vanishes and has simple poles at
the even negative integers, in particular at $s=0$. Since
$B(\pi_\infty; Ad)$ contains $0$, $\Gamma_\R(s)$ is a factor of
$L(s, \pi_\infty; Ad)$. We must then have
$$
-{\rm ord}_{s=0} L(s, \pi_\infty; Ad) \, \geq \, 1,
$$
as asserted.

\bigskip

{\it Proof of Proposition B} (contd.) \, As the local factors
$L(s, \pi_v; Ad)$ never vanish, and since $S$ is finite by
assumption, the function
$$
L_{S-\{\infty\}}(\pi; Ad) :\, = \, \prod\limits_{p \in
S-\{\infty\}} L(s, \pi_p; Ad)
$$
is non-zero at $s=0$. Hence by Lemma 2.6,
$$
-{\rm ord}_{s=0}L_S(s, \pi; Ad) \, \geq \, 1. \leqno(2.10)
$$
But we know that the full adjoint $L$-function $L(s, \pi; Ad)$ has
no pole at $s=1$, and hence at $s=0$ by the functional equation.
So all this forces the following:
$$
L^S(0, \pi; Ad) \, = \, 0,
$$
which contradicts Lemma 2.5. The only unsupported assumption we
made was that $L^S(s, \pi; Ad)$ is of positive type, which must be
wrong if $S$ is finite. We are done.

\noindent{QED}

\medskip

Note that in the proof uses the base field $\Q$ in order to use
the crucial property of the Riemann zeta function that it does not
vanish in the real interval $(0,1)$. For general number fields
$F$, the Dedekind zeta function $\zeta_F(s)$ should not have any
such real zero either, save possibly at $s=1/2$. Clearly, Theorem
A will follow for any $F$ for which Proposition B can be
established. One way to get around the difficulty for general $F$
would be to prove {\it a priori} that the adjoint $L$-function,
which has been studied from different points of view by
D.~Ginzburg, Y.~Flicker, H.~Jacquet, S.~Rallis, F.~Shahidi and
D.~Zagier, has no pole in $(0,1)$, which is, to our knowledge,
unknown. To elaborate a little, a particular version of the trace
formula due to H.~Jacquet and D.~Zagier ([JZ]) suggests that the
divisibility of $L(s, \pi \times \overline \pi)$ by $\zeta_F(s)$
for all cuspidal automorphic representations $\pi$ of GL$(n,
\A_F)$ with trivial central character is equivalent to the
divisibility of $\zeta_K(s)$ by $\zeta_F(s)$ for all commutative
cubic algebras $K$ over $F$. Since the latter is known for $n =3$,
one hopes that the former holds. This divisibility has been
investigated by relating it to an Eisenstein series on $G_2$ by
D.~Jiang and S.~Rallis ([JiR]), and the desired result could be
close to being established in the $n=3$ case.

\vskip 0.2in

\section*{\bf Bibliography}

\begin{description}

\item[{[BRa]}] L.~Barthel and D.~Ramakrishnan, A non-vanishing result for twists of
$L$-functions of GL$(n)$, Duke Math. Journal {\bf 74}, no.3,
681--700 (1994).

\item[{[C$\ell$]}] L.~Clozel, Repr\'esentations galoisiennes associ\'ees aux repr\'esentations
automorphes autoduales de ${\rm GL}(n)$, IHES Publications Math.
{\bf 73} (1991), 97--145.

\item[{[HaT]}] M.~Harris and R.~Taylor, The geometry and cohomology of some
simple Shimura varieties, with an appendix by V.~Berkovich, Annals
of Math. Studies {\bf 151}, Princeton (2001).

\item[{[He]}] G.~Henniart, Une preuve simple des conjectures de Langlands pour ${\rm GL}(n)$
sur un corps $p$-adique, Inventiones Math. {\bf 139} (2000), no.
2, 439--455.

\item[{[HRa]}] J.~Hoffstein and D.~Ramakrishnan, Siegel zeros and cusp forms,
International Math. Research Notices (IMRN) {\bf 1995}, no. 6,
279--308.

\item[{[JPSS]}] H.~Jacquet, I.~Piatetski-Shapiro and J.~Shalika,
Rankin-Selberg convolutions, Amer. J. Math.{\bf 105}, 367--464
(983).

\item[{[JS]}] H.~Jacquet and J.A.~Shalika, On Euler products and
the classification of automorphic forms I \& II, Amer. J of
Math. {\bf 103} (1981), 499--558 \& 777--815.

\item[{[JZ]}] H.~Jacquet and D.~Zagier, Eisenstein series and the
Selberg trace formula II, Transactions of the AMS {\bf 300} (1)
(1987), 1--48.

\item[{[JiR]}] D.~Jiang and S.~Rallis, Fourier coefficients of the
Eisenstein series of the exceptional group of type $G_2$, Pacific
Journal of Math. {\bf 181}, no. 2 (1997), 281--314.

\item[{[KSh]}] H.~Kim and F.~Shahidi, Cuspidality of symmetric powers with applications,
Duke Math. Journal {\bf 112} (2002), no. 1, 177--197.

\item[{[MW]}] C.~Moeglin and J.-L.~Waldspurger, Poles des fonctions $L$ de paires pour GL$(N)$,
Appendice, Ann. Sci. \'Ecole Norm. Sup. (4) {\bf 22}, 667--674
(1989).

\item[{[Pic]}] {\it Zeta functions of Picard modular surfaces},
ed. by R.P.~Langlands and D.~Ramakrishnan, CRM Publications,
Montr\'eal (1992).

\item[{[Ra]}] D.~Ramakrishnan, On the coefficients of cusp forms,
Math Research Letters {\bf 4} (1997), nos. 2--3, 295--307.

\item[{[Sh1]}] F.~Shahidi, On certain $L$-functions, American
Journal of Math. {\bf 103} (1981), 297--355.

\item[{[Sh2]}] F.~Shahidi, On the Ramanujan conjecture and the
finiteness of poles for certain $L$-functions, Ann. of Math. (2) {\bf 127}
(1988), 547--584.

\item[{[Sh3]}] F.~Shahidi, A proof of the Langlands conjecture on
Plancherel measures; Complementary series for $p$-adic groups,
Ann. of Math. {\bf 132} (1990), 273--330.

\end{description}

\vskip 0.3in

\noindent Dinakar Ramakrishnan

\noindent 253-37 Caltech, Pasadena, CA 91125

\noindent dinakar@its.caltech.edu
\end{document}